\documentclass[11pt]{article}
\usepackage[numbers]{natbib}
\usepackage{amssymb}
\begin{document}

\title{\bf Higher derivative discontinuous solutions 
to linear ordinary differential equations: A new route 
to complexity? 
}

\author{
Dhurjati Prasad Datta\thanks{Corresponding author, 
email:dp${_-}$datta@yahoo.com} and Manoj Kumar 
Bose\thanks{email:manojkbose@rediffmail.com}\\
Department of Mathematics, North Bengal University,\\
P.O. North Bengal University, 
Darjeeling, India, Pin: 734430} 
\date{}

\maketitle

\begin{abstract}
We present a new one parameter family of second derivative discontinuous 
solutions to the simplest scale invariant linear ordinary differential 
equation. We also point out how the construction could be extended to 
generate families of higher derivative discontinuous solutions as well. 
The discontinuity can occur only for a subset of even order derivatives, 
viz.,2nd, 4th, 8th, 16th, ....The solutions are shown to break the 
discrete parity (reflection) symmetry of the underlying equation. These 
results are expected to gain significance in the contemporary search of a 
new {\em dynamical principle} for understanding complex phenomena in 
Nature. 
\end{abstract}

\begin{center}
\em {Chaos, Solitons and Fractals, vol 22, issue 2, (2004), 271-275}
\end{center}

\newpage

\section{Introduction}

\par Let us consider the simplest scale invariant linear ordinary 
differential equation (ODE) 

\begin{equation}\label{ss}
t{\frac{{\rm d}\tau}{{\rm d}t}}=\tau
\end{equation}

\noindent A well known fact of the standard analysis (calculus) is that 
this equation , and more generally, an ODE of the form ${\frac{{\rm 
d}\tau}{{\rm d}t}}=f(\tau, t)$ possesses a unique $C^r$ solution in a 
neighbourhood of $(t_0,\tau_0),\,\tau_0=\tau(t_0)$ when $f\in C^r$. This 
is the Picard's theorem. The standard solution of eq(\ref{ss}), as 
admitted by the Picard's theorem, is of course $\tau_s(t)=t$ when the 
initial condition is chosen as $\tau(1)=1$. Clearly, this is a $C^\infty$ 
solution. In the following section, we, however, present, in 
contradiction  to the Picard's theorem, a new exact solution of 
eq(\ref{ss}), which is only $C^1$. The new solution breaks the reflection 
symmetry ($t\rightarrow -t$) of the underlying ODE spontaneously. We 
interprete the solution in Sec. 3 and show how this defines a nonstandard 
extension of the real number system. Consequently, every real number is 
identified with an equivalence class of a continuum of new, 
infinitesimally separated elements, which are in a state of random 
fluctuations. We also point out how the new construction generalizes to 
create higher derivative (4th, 8th,...) discontinuous families of 
solutions to the above equation. The results presented here are new and 
sums up our earlier investigations~\cite{dp1,dp2} of eq(\ref{ss}) in a 
mathematically rigorous manner. These rigorous results should become 
important in the contemporary search of a {\em new dynamical principle} 
for understanding  complexity at different levels of life viz., physical, 
biological, financial and social \cite{prigo, west,nash}, besides their 
deep significance in mathematical analysis. Ref\cite{dp2} contains some 
preliminary results and definitions indicating  their role in defining a 
new framework for addressing {\em living systems} (see also 
\cite{weiss}).

\section{$C^1$ solution}

\par Let us first demonstrate the existence of such a solution in the 
neighbourhood of $t=1$. The solution would only have a continuous  first 
order derivative at $t=1$. We begin  by defining our notations. 

\par Let $t_{n\pm}=1\pm \eta_n,\, t_0 \equiv t,\, 0<\eta_n<<1, \, 
\alpha_n=1+\epsilon_n,\, n=0,1,2,\ldots,\, {\rm and}\,\, \epsilon 
_0=0,\,0 < \epsilon_n <1 (n \neq 0) $, such that $\epsilon_n \rightarrow 
0, \, {\rm as} \, n \rightarrow \infty$ (we retain the symbol $\alpha_0$ 
for the sake of symmetry). Next, we write $t_{n\pm}^\prime=1 \pm \alpha_n 
\eta_n^\prime$, so that $\alpha_n t_{n-} =t_{n-}^\prime$. Consequently, 
$\eta_n^\prime=\eta_n - {\frac{\epsilon_n}{\alpha_n}}$. Here, $\alpha_n$ 
(and $\epsilon_n$) are scaling parameters. A useful choice,  however, is 
$\epsilon_n=\epsilon^{2^n},\, \epsilon =\epsilon _1$ (the reason will 
become clear later). As will become evident, $\eta_{n+1}=\alpha_{n}^2 
\eta_{n}^{\prime\, 2}$.

\par To construct a nontrivial solution (with the initial condition 
$\tau(1)=1$), we begin with an initial approximate solution in the small 
scale variable $\eta_0$, and then obtain recursively self-similar 
corrections over smaller and smaller scale variables $\eta_0^2,\, 
\eta_0^4,\ldots$. To this end, let  

\begin{equation}\label{ns}
\tau (t)=
\cases{
\tau_- & {\rm if} $t\lessapprox 1$\cr
\tau_+ & {\rm if} $t\gtrapprox 1$},\,\,
\tau_{-}(t_-)=(1/t_+^\prime)f_{1-}(\eta_0),\, \tau_+(t_+)=t_+
\end{equation}

\noindent be an exact solution of eq(\ref{ss}). This is obviously true 
for the right hand component $\tau_+$.  To verify the same for the 
nontrivial component $\tau_-$, we differentiate it with respect to $t_-$, 
and use the scale invariance of eq(\ref{ss}). 
Utilizing $\alpha_0 t_-=t_-^\prime$, one obtains 

\begin{equation}\label{e1}
t_-^\prime{\frac{{\rm d}\tau_{-}}{{\rm d} t_-^\prime}}=\tau_- 
({\frac{t_-^\prime}{t_+^\prime}}-t_-^\prime
{\frac{f_{1-}^\prime}{f_{1-}}})
\end{equation}

\noindent where $f_{1-}^\prime={\frac{{\rm d}f_{1-}}{{\rm 
d}\bar\eta_0}}$, $\bar\eta_0=  \alpha_0 \eta_0^\prime$. Consequently, 
eq(\ref{ns}) would be an exact solution if and only if $f_{1-}$ solves 
exactly the     self-similar equation

\begin{equation}\label{ss1}
t_{1-}{\frac{{\rm d}f_{1-}}{{\rm d} t_{1-}}}=f_{1-} 
\end{equation}

\noindent in the smaller logarithmic variable $\ln t_{1-}^{-1}$, where 
$t_{1-}=1-\alpha_0^2 \eta_0^{\prime \, 2}\equiv 1-\eta_1$. The 
self-similar replica eq(\ref{ss1}) follows from the equality 

\begin{equation}\label{id}
{\frac{t_-^\prime}{t_+^\prime}}-t_-^\prime
{\frac{f_{1-}^\prime}{f_{1-}}}=1
\end{equation}

\noindent so that $\tau_-$ is an exact solution of eq(\ref{ss}). The 
exact (nontrivial part of the ) solution could thus be written 
recursively as 

\begin{equation}\label{recu}
\tau_-=C{\frac{1}{t_+}}{\frac{1}{t_{1+}^\prime}}\ldots
{\frac{1}{t_{{(n-1)}+}^{\prime}}}f_{n-}(\eta_n^\prime)
\end{equation}

\noindent where $f_{n}$ satisfies the $n$th generation self-similar 
equation  

\begin{equation}\label{ssn}
t_{n-}{\frac{{\rm d}f_{n-}}{{\rm d} t_{n-}}}=f_{n-}
\end{equation}

\noindent and $t_{n-}=1-\alpha_{n-1}^2 \eta_{n-1}^{\prime \, 2}\equiv 
1-\eta_n$. We also note that $t_+^\prime= t_+$, since $\alpha_0=1$.

\par Plugging in the initial condition $\tau_{\pm} =1$ at $t_{\pm} =1$ 
(viz., $\eta_0 =0$), one obtains finally the desired solution as 

\begin{equation}\label{ns1}
\tau_-=C{\frac{1}{t_+}}{\frac{1}{t_{1+}^\prime}}{\frac{1}{t_{2+}^{\prime}
}}
\ldots, \, \, \, \tau_+=t_+
\end{equation}

\noindent where $C=t_{1+}^\prime(0)t_{2+}^\prime(0)\ldots$. Notice that 
$C\neq 1$, since $ \eta_1^\prime= -\epsilon_1/\alpha_1, \,  
\eta_2^\prime=\epsilon_1^2 -\epsilon_2/\alpha_2, $ etc, when $\eta_0 =0$.   

\par A remark is in order here. 
\par The solution eq(\ref{ns1}) follows from eq(\ref{recu}) only if the 
sequence $\{f_{n-}\}$ is convergent. In fact, we prove that $f_\infty = 
\lim _{n \rightarrow \infty}f_{n-} = 1$. Let $\tau_n = 
{\frac{1}{t_+}}{\frac{1}{t_{1+}^\prime}}\ldots
{\frac{1}{t_{n+}^{\prime}}}$. Then for $\eta_0$ sufficiently small and 
$\epsilon_n \rightarrow 0$ for $n \rightarrow \infty$, the sequence 
$\{\tau_n \}$ is convergent (to a nonzero value), since $t_{n+}^\prime 
\rightarrow 1$ as $n \rightarrow \infty$. Accordingly, for $\varepsilon > 
0$,  $\exists N_1$ such that $|\tau_m - \tau_n|< \varepsilon$ for $m,\, n 
> N_1 \, (m > n)$. As a result, $0< k_1 < \tau_n < k_2,\,\, k_1,\, k_2 
\sim$  O(1), for  $n > N_2$ for a sufficiently large $N_2$. Again, 
$f_{n-}$, being defined by eq(\ref{ssn}), is uniformly bounded in a 
neighbourhood of $t=1$, so that $|f_{n-}|< k$ for $n > N_2$. The desired 
convergence now follows from the Cauchy convergence criterion, since 
$|f_{n-} - f_{m-}| = |\tau_n^{-1}||f_{m-}||\tau_m - \tau_n| < k_1^{-1} k 
\varepsilon$ $ \forall m, n > N,\, N$= max($N_1,N_2$). Finally, 
eq(\ref{ssn}), in the asymptotic limit $n \rightarrow \infty$, yields 
$f_\infty = {\frac{{\rm d}f_{n-}}{{\rm d} t_{n-}}}|_\infty ={\frac{{\rm 
d}\tau}{{\rm d} t}}|_{t=1} = \tau(1)=1$.      

\par Now to test the continuity of the derivatives of the solution 
(\ref{ns1}) at $t_\pm=1$, i.e., at $\eta_0=0$, we note that 
$\eta_n^\prime$ is a polynomial in $\eta_0$, of degree $2^n$, being 
defined recursively by $\eta_n^\prime=\eta_n - 
{\frac{\epsilon_n}{\alpha_n}},\, \eta_n=\alpha_{n-1}^2 \eta_{n-1}^{\prime 
\,2}$. As a result ${\frac{{\rm d}\eta_n^\prime}{{\rm d}\eta_0}}=0$, but 
${\frac{{\rm d^2}\eta_n^\prime}{{\rm d}\eta_0^{ 2}}} \neq 0$, at 
$\eta_0=0$. One thus obtains

\begin{equation}\label{d1}
{\frac{{\rm d}\tau_-}{{\rm d}\eta_0}}=- \tau_- \{{\frac{1}{1+ \eta_0}} + 
({\frac{\alpha_1}{1+ \alpha_1 \eta_1^\prime}}) {\frac{{\rm 
d}\eta_1^\prime}{{\rm d}\eta_0}} + ({\frac{\alpha_2}{1+ 
\alpha_2\eta_2^\prime}}) {\frac{{\rm d}\eta_2^\prime}{{\rm d}\eta_0}} + 
\ldots \}
\end{equation}  

\noindent so that ${\frac{{\rm d}\tau_-}{{\rm d}t_-}}=1={\frac{{\rm 
d}\tau_+}{{\rm d}t_+}}$ at $\eta_0 =0$ which means that the first 
derivative of the solution is indeed continuous for all $\eta_0$. 
However, as is verified easily from eq(\ref{d1}), the second derivative 
of $\tau_-$ at $\eta_0 =0$ is not zero, as one expects on the basis of 
the standard solution $\tau_s =t$. Indeed, one can verify that 
${\frac{{\rm d^2}\tau_-}{{\rm d}\eta_0^{2}}}= 
2(1-{\frac{1+\epsilon_1}{1-\epsilon_1}} - \ldots) \neq 0$ at $\eta_0 =0$, 
unless $\epsilon_n =0$, for all $n$. In this special case, i.e., when 
$\epsilon_n=0,\, \forall n$, our solution (\ref{ns1}) reduces to the 
standard solution, since $\tau 
_-={\frac{1}{1+\eta_0}}{\frac{1}{1+\eta_0^2}}{\frac{1}{1+\eta_0^4}}\ldots 
=1-\eta_0=t_-$.

\par It thus follows that the solution (\ref{ns1}), with nonzero scaling 
parameters, is indeed nontrivial, because of this second derivative 
discontinuity at $\eta_0 =0$, that is at $t =1$. In fact, the scaling 
invariance of eq(\ref{ss}) tells also that, $t=1$ could be realized as $t 
\rightarrow t/t_0=1$, so that the nontrivial solution (\ref{ns1}) 
actually holds in the neighbourhood of every real number $t_0$, the 2nd 
derivative being discontinuous at $t=t_0$. This is our main result of 
this paper. Let us note here that $\tau_-=\tau_{s-}(1 + {\rm 
O}(\eta_0^2))$, besides the arbitrariness of the scaling parameters 
$\epsilon_n$. Combining the standard and the new solutions together, one 
can finally write down a more general one parameter class of solutions 

\begin{equation}\label{gen}
\tau_g(t)=t(1 + \phi(t)), \, \, \phi(t)=\epsilon t^{-1} \tau(t)
\end{equation} 

\noindent Note that  

\begin{equation}\label{gc}
t{\frac{{\rm d}\phi}{{\rm d}t}}=0
\end{equation}

\noindent because $\tau$ is an exact solution of eq(\ref{ss}) (we fix 
$\epsilon_n = \epsilon^{2^n}$, on anticipation of Sec.3). The 2nd 
derivative discontinuity of $\tau$, however, tells that $\phi$ can not be 
considered simply as an ordinary constant~\cite{dp1, dp2}.  
\subsection{Symmetry breaking}

\par The solution (\ref{ns1}) is not parity (reflection) symmetric. Let 
$P: P t_\pm =t_\mp$ denote the reflection transformation near $t=1$ 
($P\eta=-\eta$ near $\eta=0$). Clearly, eq(\ref{ss}) is parity symmetric. 
So is the standard solution $\tau_{s\pm} = t_\pm$ (since $P \tau_s = 
\tau_s$). However, the solution (\ref{ns1}) breaks this discrete symmetry 
spontaneously: $\tau_-^P=P \tau_+ = t_-,\, \tau_+^P = P \tau_- = 
C{\frac{1}{t_-}} {\frac{1}{t_{1+}^\prime}} 
{\frac{1}{t_{2+}^{\prime}}}\ldots$, which is of course a solution of 
eq(\ref{ss}), but clearly differs from the original solution, $\tau_\pm^P 
\neq \tau_\pm$. The explicit breaking of the parity symmetry is an 
important new feature of the solution (\ref{ns1}). Its significance in 
dynamical systems would be studied elsewhere. 

\section{Interpretation}

\par A basic assumption in the framework of the standard calculus is that 
a real variable $t$ changes by linear translation only. Further, $t$ 
assumes (attains) every real number exactly. However, in every 
computational problem within a well specified error bar, a real number is 
determined only upto a finite degree of accuracy $\epsilon_0$ ,say. 
Suppose,  for example, in a computation, 1 is determined upto an accuracy 
of $\pm 0.01$, so that 1  here effectively   stands for the set 
$1_\epsilon\equiv\{1 \pm \epsilon,\,\epsilon < \epsilon_0 = 0.01\}$ with 
cardinality $c$, the cardinality of the continuum. We call $\epsilon$ an 
`infinitesimally small' real number. Now, any laboratory computational 
problem (experiment) is run only over a finite time span, and the 
influence of such infinitesimally small $\epsilon$'s, being 
insignificantly small, could in fact be disregarded. However, if the 
experiment is allowed to run over an `infinitely' large period of time, 
analogous to a natural (or biological/ financial) process, such as the 
fluctuations of a river (ocean) water level (heart beats of a mammal/ 
stock market prices ), which is known to occur over many time scales, the 
contributions from the infinitesimally small scales $\epsilon$ need not 
remain  negligible. 

\par Our solution (\ref{ns1}) captures the essence of the above scenario 
in an elegant manner. Notice that our initial ansatz $\tau_-\approx 
1/t_+$, approximates the standard solution $\tau_{s-}=t_-$ upto 
O($\eta_0^2)$ (i.e., when $\eta_0\neq 0$, but (O$(\eta_0^2)=0$). At this 
level of approximation $1_\epsilon=1$ for $\epsilon \lessapprox \eta^2$. 
However, at the second level of iteration with the variable 
$t_{1-}=1-\eta_0^2$, the effect of $\epsilon$, being of the same order of 
$\eta_0^2$, can not be negligible. Recalling that 1 here should be 
identified as $1 \pm \epsilon$, an appropriate variable can therefore be 
written as  $\alpha t_{1-}=1 \pm \epsilon -\alpha\eta_0^2$, which in turn 
is  written equivalently as the rescaled variable  $t_{1-}^\prime = 
1-\alpha \eta_1^\prime$ where $\eta _{1} ^\prime =\eta_0^2 \mp 
\epsilon/\alpha $ so that $t_{1-}^\prime = \alpha t_{1-}, \, \alpha = 1 
\pm \epsilon$.
This scaling transformation, which is used recursively in (\ref{ns1}), 
thus not only reveals an intrinsically  approximate nature of 
$1_\epsilon$, but also provides an window to probe the fine scale 
structure of $1_\epsilon$ by recursively approaching more and more 
accurately the `exact' value 1 through the finer and finer scales 
$\epsilon^2,\, \epsilon^4, \ldots$ successively. ( This explains the 
choice $\epsilon_n = \epsilon^{2^n}$.) The exact value of 1 (and, hence, 
of any real number), as endowed in the framework of standard classical 
analysis, is, however, {\em illusory}, because of the fact that the 
solution (\ref{ns1}), being 2nd derivative discontinuous at $t=1$, is a 
{\em new exact} solution of eq(\ref{ss}). Further, the value of a small 
real number $\eta_0$ is uncertain not only upto O($\eta_0^2$), but also 
because of the arbitrary parameter $\epsilon$. We note that two solutions 
$\tau_g$ and $\tau_s$ are indistinguishable for $t \sim$ O(1) and 
$\eta_0^2 <<1$. However, for a sufficiently large $t \sim$ 
O($\epsilon^{-1})(\equiv$ O($\eta_0^{-2}$)), the behaviours of two 
solutions would clearly be different. Moreover, because of  the intrinsic 
(irreducible) uncertainty {\em the new solution $\tau_g$ could 
effectively be written as $\tau_g:\tau_-\approx 1/t_+,\,\tau_+=t_+$, 
without any loss of generality (and information) at the order $\eta_0\neq 
0,\,\eta_0^2=0$}. The correction factors in eq(\ref{ns1}) only reveals 
the approximate (statistical) self similarity of the solution over higher 
order scales $\eta_0^2,\,\eta_0^4,\ldots$ and so on.    

\par To summarise, the nature of the solution (\ref{ns1}) depends 
crucially on the nontrivial scaling $\alpha_n t_-=t_-^\prime$, as defined 
above, which shifts $\eta_0=0$ at the level of first iteration to 
$\eta_1^\prime =0$ for the second iteration and so on recursively to 
finer and finer scales. The origin of the second derivative discontinuity 
is obviously an effect of this receding of the zero (0) down finer 
scales. Let us also point out that the scaling  $\alpha t_-=t_-^\prime$, 
does not mean $\alpha t_+=t_+^\prime$ (in fact, $t_+^\prime = 
(1-\epsilon) (1 + {\frac {1 + \epsilon} {1 - \epsilon}} \eta_0)\neq 
\alpha t_+$), which in turn leads to the nontrivial feature that is 
revealed in the solution (\ref{ns1}).

\par Now to state the above more formally, we note that the standard 
solution $\tau_s=t$ defines the 1-1 identity map  $R \rightarrow R$ of 
the real number set $R$. The generalized solution (\ref{gen}) now 
provides a one parameter extension of the identity map ${\bf 1}_\epsilon: 
R \rightarrow R$. Accordingly, every real number $t$ is mapped to the 
corresponding `fat' real number ${\bf t}_\epsilon= t (1+\phi)$, which  is 
the set of `infinitely many' {\em new} elements (with cardinality $2^c$) 
separated by infinitesimally small numbers (distances) \cite{dp1}, since  
$\phi$  is infinitesimally small. Accordingly, the set of real numbers 
$R$ is considered to be an enlarged nonstandard set \cite{robin} ${\bf 
R}$ with infinitesimally small elements, which remain dormant in the 
standard analysis (i.e., in the context of the standard solution $\tau_s$ 
of eq(\ref{ss})). This extension of $R$ to the nonstandard set ${\bf R}$ 
via the generalized solution of the linear ODE eq(\ref{ss}), is analogous 
to the extension of the ring of integers to the field of rational numbers 
as the solution space of the linear equation $ax + b=0,\, a,\, b$ being 
integers. Let us note further that eq(\ref{d1}) written as 

\begin{equation}\label{d2}
{\rm d}\tau_- = - \tau_- \{{\frac{{\rm d} \eta_0}{1+ \eta_0}} + 
{\frac{\alpha_1 {\rm d}\eta_1^\prime}{1+ \alpha_1 \eta_1^\prime}} + 
{\frac{\alpha_2 {\rm d}\eta_2^\prime}{1+ \alpha_2\eta_2^\prime}} + \ldots 
\}
\end{equation}  

\noindent tells that the variables $ t_{+},\, t_{1+}^\prime, \ldots$ 
behave as independent variables. Writing ${\bf t}_-^{-1} = \Pi_0^{\infty} 
t_{n +}^\prime$, a fat real variable in the neighbourhood of ${\bf 
1}_\epsilon$, eq(\ref{d2}) restates the fact that the solution 
eq(\ref{ns1}) is indeed an `exact' solution of the equation $ {\frac{{\rm 
d}\tau}{{\rm d}\ln {\bf t_-}}}=\tau$ in ${\bf R}$. In $R$, however, we 
have an irreducible `unknown' component $f_{n-}\sim$O(1), which is 
exactly determined only in the limit $n \rightarrow \infty$. Ironically, 
however, even in this asymptotic limit, one fails to improve (determine) 
the solution with an unlimited accuracy, because  of the indeterminate 
parameter $\epsilon$, as reflected in the 2nd derivative discontinuity. 
The fat variable reduces to the ordinary real variable $t_-$ only in the 
approximation when $\eta_0^2 = 0 \,(\eta_0\neq 0) $.         
 
 \par The solution (\ref{ns1}) also tells that a  fat real variable ${\bf 
t}=t{\bf 1}$ in ${\bf R}$ can change not only by linear translations but 
also by inversions: $t_-\rightarrow t_-^{-1} =t_+$, in the neighbourhood 
of $t=1$. Clearly, the inversion is realized `exactly' in {\bf R}, when 
an (first order) infinitesimal $\phi$ is defined by $\phi \neq 0$, but 
O($\phi^2$)=0. Consequently, the change (increment) of a real variable 
could be visualized as an SL(2,R) group action, viz., a combination of 
linear translation ( for ordinary real variable $t$ following  the 
standard solution of eq(\ref{ss})) and inversion ( for $t_\pm \in {\bf 
1}$). As note already, because of the arbitrariness in $\epsilon$, as 
well as being defined upto O($\phi^2=0$), the infinitesimal elements 
$\phi$ enjoy an element of randomness, which in turn renders ${\bf t}$ 
random. The ordinary non-random real variable $t$ is thus retrieved only 
under an approximation, that is, when $\phi=0$. 

\par Let us remark finally that one can generate 4th or higher derivative 
discontinuous solutions of eq(\ref{ss}) by introducing infinitesimal 
scalings at an appropriate level of iteration. Note that the 2nd 
derivative discontinuous solution is obtained when the nontrivial scaling 
viz., $t_{1-} \rightarrow t_{1-}^\prime$ is introduced at the 1st 
iteration (the ansatz in eq(\ref{ns}) stands for the zeroth iteration ). 
Thus, the 4th derivative discontinuous solution is obtained when the 
rescaling of the appropriate variable is used at the 2nd level of the 
iteration, viz., $t_{2-} \rightarrow t_{2-}^\prime$, giving rise to yet 
another class of solutions 
\begin{equation}\label{recu1}
\tau_-^{(2)}=C{\frac{1}{t_+}}{\frac{1}{t_{1+}}}{\frac{1}{t_{2+}^\prime}}\%
ldots
{\frac{1}{t_{{(n-1)}+}^{\prime}}}f_{n-}(\eta_n^\prime)
\end{equation}

\noindent where $t_{1+}= 1+\eta_0^2,\, t_{2+}^\prime = 1 - \epsilon + 
(1+\epsilon) \eta_0^4, \, t_{3+}^\prime = 1-\epsilon^2 + (1 
+\epsilon^2)[(1 + \epsilon)\eta_0^4 - \epsilon]^2,\ldots$, and so on for 
higher derivative discontinuous solutions. Note that $\tau_- 
^{(2)}=\tau_s (1 + $O($\eta_0^4$)) for an $0 \neq \epsilon \lessapprox$ 
O($ \eta_0^4$). Consequently, $\ln (\tau_-^{(n)}/\tau_s)\approx 
\eta_0^{2^{n}} $ corresponds to the $n$th order infinitesimals in {\bf 
R}.

\section{Conclusion} 
\par  We have presented  new families of higher derivative discontinuous 
solutions of  the ODE (\ref{ss}), which do not respect the Picard's 
theorem.  The origin of this violation could be traced to the fact that a 
variable in ${\bf R}$ may undergo changes (increments) via the extended 
SL(2, R)-like group actions. These solutions break explicitly the parity 
symmetry of the underlying ODE. The derivation of such solutions from  
any linear first order ODE is obvious. Higher order equations will be 
considered elsewhere. We close with the observation that the results 
presented here may as well be considered to offer, so to speak, an exact 
proof of the fact that mathematics is inexact.

\end{document}